\documentclass[10pt]{amsart}
\usepackage{latexsym, amsmath}

\renewcommand{\epsilon}{\varepsilon}
\newcommand{\newsection}[1]
{\subsection{#1}\setcounter{theorem}{0} \setcounter{equation}{0}
\par\noindent}

\newtheorem{theorem}{Theorem}

\newtheorem{lemma}[theorem]{Lemma}
\newtheorem{corr}[theorem]{Corollary}

\newtheorem{proposition}[theorem]{Proposition}
\newtheorem{deff}[theorem]{Definition}

\newcommand{\bth}{\begin{theorem}}
\newcommand{\ble}{\begin{lemma}}
\newcommand{\bcor}{\begin{corr}}

\newcommand{\bdeff}{\begin{deff}}

\newcommand{\bprop}{\begin{proposition}}
\newcommand{\eth}{\end{theorem}}
\newcommand{\ele}{\end{lemma}}
\newcommand{\ecor}{\end{corr}}
\newcommand{\edeff}{\end{deff}}

\newcommand{\eprop}{\end{proposition}}

\newcommand{\cd}{\, \cdot\, }

\newcommand{\Rn}{{\mathbb R}^n}

\renewcommand{\Pi}{\varPi}

\renewcommand{\epsilon}{\varepsilon}

\newcommand{\MU}{\mu}

\newcommand{\R}{{\mathbb R}}

\begin{document}

\title[Nonlinear equations in waveguides]
{Nonlinear hyperbolic equations in infinite homogeneous
waveguides}
\thanks{The authors were supported in part by the NSF}

\author{Jason Metcalfe}
\address{School of Mathematics, Georgia Institute of Technology, Atlanta, GA 30332-0160}
\author{Christopher D. Sogge}
\address{Department of Mathematics,  Johns Hopkins University,
Baltimore, MD 21218}
\author{Ann Stewart}
\address{Department of Mathematics,  Johns Hopkins University,
Baltimore, MD 21218}

%\maketitle

\begin{abstract}

In this paper we prove global and almost global existence theorems
for nonlinear wave equations with quadratic nonlinearities in
infinite homogeneous waveguides.  We can handle both the case of
Dirichlet boundary conditions and Neumann boundary conditions.  In
the case of Neumann boundary conditions we need to assume a
natural nonlinear Neumann condition on the quasilinear terms.  The
results that we obtain are sharp in terms of the assumptions on
the dimensions for the global existence results and in terms of
the lifespan for the almost global results.  For nonlinear wave
equations, in the case where the infinite part of the waveguide has
spatial dimension three, the hypotheses in the theorem concern
whether or not the Laplacian for the compact base of the waveguide has a zero mode or not.

\end{abstract}
\maketitle

\newsection{Introduction}

In this paper we shall consider nonlinear wave and Klein-Gordon
equations of the form
\begin{equation}\label{1.1}\begin{cases}
(\square +m^2)u=Q(u,u',u''),
\\
u(0,x)=f(x), \quad \partial_tu(0,x)=g(x).
\end{cases}
\end{equation}
where
\begin{equation}\label{1.2}
\square = \partial_t^2 - (\Delta + \Delta_\Omega)
\end{equation}
is the d'Alembertian on $(t,x,y)\in {\mathbb R}_+ \times {\mathbb
R}^n \times \Omega$.  Here,
$$\Delta = \Delta_{\Rn}
=\sum_{j=1}^n \partial^2/\partial x_j^2$$ is the Euclidean
Laplacian on ${\mathbb R}^n$. Also, $\Omega\subset {\mathbb R}^d$
denotes a nonempty, bounded domain with smooth boundary $\partial
\Omega$, and $\Delta_\Omega$ denotes the standard Laplacian
$$\Delta_\Omega=\sum_{j=1}^d \partial^2/\partial y_j^2.$$
So we shall impose either Dirichlet boundary conditions
%
% while
%$\Delta_\Omega$ is the Laplace-Beltrami associated with a smooth
%Riemannian metric $g$ on our fixed connected compact manifold $M$
%of dimension $d$ that may or may not have boundary. Also, we shall
%be able to prove natural results for the Klein-Gordon and the wave
%equation case where, in \eqref{1.1}, $m>0$ and $m=0$,
%respectively.
%
%If $\partial \Omega\ne \emptyset$ then we shall assume that the
%equation \eqref{1.1} has either Dirichlet boundary conditions
\begin{equation}\label{1.3} u(t,x,y)|_{y\in\partial \Omega}=0
\end{equation}
or Neumann boundary conditions
\begin{equation}\label{1.4}
\partial_\nu u(t,x,y)|_{y\in\partial \Omega}=0,
\end{equation}
with $\partial_\nu$ denoting the normal derivative at $\partial
\Omega$.  Thus, depending on which case we are in, $\Delta_\Omega$
will be either the Dirichlet Laplacian or the Neumann Laplacian on
$\Omega$.

Also, in \eqref{1.1}, $Q$ is a quadratic function in its
arguments, which are $u$, $u'=\nabla u$, and $u''=\nabla^2 u$,
with $\nabla=(\partial_t,\partial_x, \partial_y)$ being the full
space-time gradient.  We shall  assume that our nonlinear
equations are quasilinear, which means that they are affine linear
in $u''$. For simplicity we shall just state things for the scalar
case, but everything carries over without any difficult to the
case where \eqref{1.1} is a nonlinear system in $u=(u^1,\dots,
u^N)$.  So, we are assuming that $Q$ can be written as
\begin{equation}\label{1.4'}
Q(u,u',u'')=\sum_{0\le j,k,l\le n+d}A^{jk}_l \partial_l u
\partial_j\partial_k u
+ u\sum_{0\le j,k\le n+d}A^{jk}\partial_j\partial_ku
+R(u,u')
\end{equation}
where $R$ is a quadratic function in $u$ and $u'$.
Here, and in what follows, $x_0=t$, and $x_{n+j}=y_j, \, 1\le j\le
d.$

In order to solve \eqref{1.1} we must assume that the data
satisfies the relevant compatibility conditions. Since these are
well known (see, e.g., \cite{KSS}), we shall describe them briefly.
To do so we first let $J_ku =\{\partial^\alpha_xu: \, 0\le
|\alpha|\le k\}$ denote the collection of all spatial derivatives
of $u$ of order up to $k$, using local coordinates in a small
tubular neighborhood of $\partial \Omega$. Then if $N$ is fixed
and if $u$ is a formal $H^N$ solution of \eqref{1.1} we can write
$\partial_t^ku(0,\cdot)=\psi_k(J_kf,J_{k-1}g)$, $0\le k\le N$, for
certain compatibility functions $\psi_k$ which depend on the
nonlinear term $Q$ as well as $J_kf$ and $J_{k-1}g$. Having done
this, the compatibility condition for the Dirchlet case
\eqref{1.1}, \eqref{1.3} with $(f,g)\in H^N\times H^{N-1}$ is just
the requirement that the $\psi_k$ vanish on ${\mathbb R}_+\times
{\mathbb R}^n\times \partial\Omega$
%when $y\in \partial \Omega$
when $0\le k\le N-1$. Additionally, we shall say that
$(f,g)\in C^\infty$ satisfy the compatibility conditions to
infinite order if this condition holds for all $N$.  For the
Neumann case \eqref{1.1}, \eqref{1.4} one follows a similar
construction by writing $\partial_t^k\partial_\nu
u(0,\cd)=\Psi_k(J_{k+1}f, J_kg)$ and requiring that the $\Psi_k$
vanish on ${\mathbb R}_+\times
{\mathbb R}^n\times \partial\Omega$
% when $y\in \partial \Omega$
when $0\le k\le N-2$.  %If
%$\partial \Omega=\emptyset$ then no compatability conditions are
%needed.

For simplicity, we shall also assume that the initial data has
compact support.  So we shall assume that there is a fixed
constant $B>0$ so that
\begin{equation}\label{1.5}
f(x,y)=g(x,y)=0, \quad |x|>B,
\end{equation}
where we recall that $x$ denotes the ${\mathbb R}^n$ variable.
%We shall indicate
%in remarks the situations where our results easily lead to such
%results.
In addition to \eqref{1.5} we shall also have to assume
that the initial data is small,
\begin{equation}\label{1.6}
\|f\|_{H^N({\mathbb R}^n\times \Omega)}+\|g\|_{H^{N-1}({\mathbb
R}^n\times \Omega)}\le \varepsilon,
\end{equation}
where $\|f\|_{H^N}=\sum_{|\alpha|\le
N}\|\partial^\alpha_{x,y}f\|_{L^2(\Rn\times\Omega)}$.
This equation can be relaxed in certain situations to a condition
that certain weighted-Sobolev norms are small.

%Let us now state our main results.  The first says that we always
%have small amplitude global existence for Klein-Gordon equations
%with quadratic nonlinearities when $n\ge 3$.
%
%\begin{theorem}\label{theorem1.1}  Suppose that $m>0$ in
%\eqref{1.1}, and assume that $n\ge 3$.   Assume also that the
%Cauchy data $(f,g)\in C^\infty({\mathbb R}^n\times \Omega)$
%satisfies \eqref{1.5} and \eqref{1.6} as well as the appropriate
%infinite order compatibility conditions.   It then follows that
%the corresponding nonlinear Klein-Gordon equation \eqref{1.1} with
%boundary conditions \eqref{1.3} or \eqref{1.4} has a global smooth
%solution if $N$ in \eqref{1.6} is a sufficiently large fixed
%integer and if $0<\varepsilon<\varepsilon_0$ is sufficiently
%small.
%\end{theorem}

Let us now state our main results.  The first says that for
Dirichlet boundary conditions we always have small amplitude
global existence for Klein-Gordon or wave equations with quadratic
nonlinearities when $n\ge 3$.

\begin{theorem}\label{theorem1.1}  Suppose that $m\ge0$ in
\eqref{1.1}, and assume that $n\ge 3$.   Assume also that the
Cauchy data $(f,g)\in C^\infty({\mathbb R}^n\times \Omega)$
satisfies \eqref{1.5} and \eqref{1.6} as well as the appropriate
infinite order compatibility conditions for the Dirichlet boundary
conditions \eqref{1.3}.  It then follows that the corresponding
nonlinear hyperbolic equation \eqref{1.1}, \eqref{1.3} has a global smooth
solution if $N$ in \eqref{1.6} is a sufficiently large fixed
integer and if $0<\varepsilon<\varepsilon_0$ is sufficiently
small.
\end{theorem}

The case of Neumann boundary conditions is more delicate for
several reasons.  First, to be able to have energy estimates, one
needs a natural nonlinear Neumann condition on the quasilinear
quadratic forms (see \S 5)
\begin{equation}\label{1.8'}
\sum_{0\le j,k,l\le n+d}A^{jk}_l \xi_l \eta_j\theta_k =0  \quad
\text{and} \sum_{0\le j,k\le n+d}A^{jk}\xi_j\theta_k=0,\quad  \text{if
} \, (\theta,\xi,\eta)\in X,
\end{equation}
where
$$X=\{\, (\theta,\xi,\eta): \,
\theta=(0,\dots,0,\nu_1(y),\dots,\nu_d(y)), \, \xi\cdot\theta=0,
\, \eta\cdot\theta=0, \, y\in \partial\Omega\}.$$ Thus, $\theta$
is assumed to be normal to ${\mathbb R}^{1+n}\times
\partial\Omega$, and $\xi$ and $\eta$ are both assumed to be
orthogonal to $\theta$. Note that this condition automatically
holds if the quasilinear terms only involve $\partial_j\partial_k
u$, $0\le j,k\le n$, i.e., $A^{jk}_l=0$ and $A^{jk}=0$ if $n+1\le
j\le n+d$ or $n+1\le k\le n+d$.  This condition is the natural one so
that energy estimates can hold.  It should be compared to the symmetry condition
for multispeed nonlinear hyperbolic systems (see, e.g., \cite{KSS3}, \cite{MS}, \cite{Si}).

If we assume this nonlinear condition, then in the case of
Neumann boundary conditions we get the same sort of results for
Klein-Gordon (i.e., $m>0$) equations:

\begin{theorem}\label{theorem1.2}  Suppose that $m>0$ in
\eqref{1.1}, and assume that $n\ge 3$.   Assume also that the
Cauchy data $(f,g)\in C^\infty({\mathbb R}^n\times \Omega)$
satisfies \eqref{1.5} and \eqref{1.6} as well as the appropriate
infinite order compatibility conditions for the Neumann boundary
conditions \eqref{1.4}.  Then, if the nonlinear Neumann condition
\eqref{1.8'} is satisfied, it follows that the corresponding
nonlinear Klein-Gordon equation \eqref{1.1}, \eqref{1.4} has a global smooth
solution if $N$ in \eqref{1.6} is a sufficiently large fixed
integer and if $0<\varepsilon<\varepsilon_0$ is sufficiently
small.
\end{theorem}

Proving existence theorems for nonlinear wave equations  (i.e.,
$m=0$) in waveguides with Neumann boundary conditions is more
difficult. However, we can obtain optimal results for certain
semilinear equations without too much difficulty using ideas from
\cite{KSS2}.

\begin{theorem}\label{theorem1.3} Consider the nonlinear wave
equation \eqref{1.1} with $m=0$.  Assume that the quadratic
nonlinearity is independent of $u$, $\partial_tu$, $\nabla_yu$ and
$u''$, i.e., that \eqref{1.1} is replaced by
\begin{equation}\label{1.7}\square u = Q(\nabla_xu),
\end{equation}
where $Q(\nabla_xu)$ is a bilinear form in the spatial ${\mathbb
R}^n$ gradient.
% ${\mathbb R}\times \Rn$ gradient.
 Assume that the
data $(f,g)$ are as in Theorem \ref{theorem1.2}. Then if $n\ge 4$
\eqref{1.7}, \eqref{1.4} has a global smooth solution if
$\varepsilon<\varepsilon_0$ is sufficiently small.  If $n=3$
 and $\varepsilon<\varepsilon_0$ is
sufficiently small, there is a constant $c>0$ so that if
\begin{equation}\label{1.8}
T_\varepsilon = \exp(c/\varepsilon),
\end{equation}
then \eqref{1.7}, \eqref{1.4} has a solution $u\in
C^\infty([0,T_\varepsilon) \times {\mathbb R}^3\times \Omega)$.
\end{theorem}

This result is of course optimal.  If one takes data
$f(x,y)=f(x)$, $g(x,y)=g(x)$ that are independent of $y\in \Omega$
then the solution of the equation \eqref{1.7}, \eqref{1.4} is
given by $u(t,x,y)=u(t,x)$ where $u$ is the solution of the
equation $\square u(t,x)=Q(\nabla_x u(t,x))$ in Minkowski space.
Classic results of John (see, e.g., \cite{john1,john2}) show that
there is blowup for times of order $\exp(C/\varepsilon)$
 when
$n=3$ and $Q(\nabla_xu)=|\nabla_x u|^2$.

It is for technical reasons that we can only handle semilinear
terms involving $\nabla_x u$ in dimensions $n=3,4$.  We remark,
though, that when $n\ge5$ it seems that our extension \eqref{2.6}
of the KSS inequality \cite{KSS2} does give global existence for
equations involving bilinear forms in $(\partial_tu,\nabla_x u)$.
We shall study this in a future paper.

We remark that we could also have obtained the analog of Theorem
\ref{theorem1.2} or Theorem \ref{theorem1.3} in the case where
$\Omega$ is a compact Riemannian manifold without boundary.  As we
shall see in the proof the case of Dirichlet boundary conditions
is much more favorable than the boundaryless case or the case of
Neumann boundary conditions when $m=0$.  This is because the
latter cases have zero eigenmodes.  In the case of Dirichlet
boundary conditions, since there are no zero eigenvalues, one can
naturally reduce matters to proving uniform estimates for
Klein-Gordon equations $\square_{{\mathbb R}^{1+n}} + \mu^2$ in ${\mathbb R}\times
{\mathbb R}^n$ with $\mu\ge 1$.  On the other hand, if $m=0$, then
for the Neumann or boundaryless case, one must also consider estimates
for $\square_{{\mathbb R}^{1+n}}$, which are less favorable in terms
of time-decay.  For example, see the Klainerman-Sobolev inequalities \cite{K} (see also, e.g., \cite{H} and \cite{S}) which
roughly provide $t^{-(n-1)/2}$ decay for the wave equation as compared to the $t^{-n/2}$ decay
in Proposition \ref{prop2.3} and Proposition \ref{prop2.4} for the Klein-Gordon equation.
Here, and throughout, $\square_{{\mathbb R}^{1+n}} =
\partial_t^2-\Delta_{{\mathbb R}^n}$.

In order to show global existence of solutions to the nonlinear problem \eqref{1.1} in
the proofs of Theorem \ref{theorem1.1} and Theorem \ref{theorem1.2}, we will use the above
a priori estimates and an energy estimate
in a fashion similar to that used by Klainerman \cite{kkg} in the boundaryless case.
See also H\"ormander \cite{H}.  It should be noted that Shatah \cite{shatah} provides
an alternate proof in the boundaryless case that uses normal forms to reduce to the case of
nonlinearities vanishing to third order.  From here, the proof proceeds using estimates that
are relatively straightforward as compared to the a priori estimates used in \cite{kkg}.  We, however,
do not explore Shatah's method in this paper.  As mentioned previously, the proof of global existence
for wave equations in waveguides with Neumann boundary conditions is more difficult.  A proof of Theorem
\ref{theorem1.3} which uses techniques reminiscent of those of \cite{KSS2} is provided.

Earlier results were obtained by Lesky and Racke \cite{lesky}, and
their work inspired this one.
%Using techniques of Shibata and
%Tsutsumi \cite{tsutsumi} they were able to show that the above
%global existence results hold in the case of Dirichlet boundary
%conditions when $n\ge 5$.
Their techniques largely relied on
$L^p\to L^{p/(p-1)}$ decay estimates for the linear Klein-Gordon
equation that are reminiscent to ones obtained earlier by
Marshall, Strauss and Wainger \cite{msw}.  As with our estimates,
the main point is to prove estimates for solutions of Klein-Gordon
equations $(\square_{{\mathbb R}^{1+n}} +m^2)u=0$ in ${\mathbb R}\times {\mathbb R}^n$
that are uniform in $m\ge 1$.  By using an eigenfunction expansion
for the other part of the waveguide, $\Omega$, one can use these
to easily obtain estimates on ${\mathbb R}\times ({\mathbb
R}^n\times \Omega)$.
%that are sufficient for proving the above
%existence results.
  Then, using techniques of Shibata and Tsutsumi \cite{tsutsumi}
they are able to use these linear estimates to show that the above
global existence results hold in the case of Dirichlet boundary conditions
when $n\ge 5$.
If $n$ is fixed then as the dimension of
$\Omega$ increases, one of course needs more and more regularity of the
initial data.

Our techniques are related to recent work on nonlinear obstacle
problems, (e.g., \cite{KSS}, \cite{KSS2}, \cite{KSS3}, \cite{MS},
and \cite{SS}).  However, the results here are considerably easier
to obtain for a couple of reasons.  First, and most important,
when we use the method of commuting vector fields (see, e.g., \cite{S}) in our case,
since there is no obstacle in ${\mathbb R}^n$, we are allowed to
use the generators of hyperbolic rotations, i.e.,
$$\Omega_{0j}=x_j\partial_t + t\partial_j, \quad 1\le j\le n,$$
as well as the generators of spatial rotations,
$$\Omega_{jk}=x_j\partial_k-x_k\partial_j, \quad 1\le j<k\le n.$$
 Also since there
is no obstacle in ${\mathbb R}^n$, we do not have to use the unit
cutoff method that relies on local exponential decay of energy
(\cite{Ikawa1}, \cite{Ikawa2}, \cite{LMP}, \cite{MRS}). If one
wishes to extend our results to the case of inhomogeneous wave
guides where one considers a compact perturbation of the ${\mathbb
R}^n$-Laplacian, either from a change of metric or an obstacle,
then the situation would be much more delicate.  We should also
point out that in our situation we cannot use the scaling vector
field $t\partial_t + r\partial_r$ since it does not preserve the
equation $(\partial_t^2-\Delta + \mu^2)u=0$, $\mu\ne 0$.  However,
this presents no problems since the Sobolev estimates and weighted
$L^2(dtdx)$ estimates are strong enough. The latter follow easily
from a slight variant of a local estimate of Smith and the second
author \cite{SS}.

The paper is organized as follows.  In the next two sections, we
shall state the estimates for the linearized equation that we
shall need. We first prove estimates in ${\mathbb R}^n$, including
ones which involve bounds that are uniform in the mass parameter
$m$ for the linear equation.  We then show how these uniform
bounds easily lead to estimates for solutions of linear equations
in waveguides.  Using these uniform estimates and orthogonality
we can easily obtain the existence theorems using standard
arguments.  These proofs are carried out in \S 4 for the Dirichlet
case and in \S 5 for the Neumann cases.
In what follows, $C$ shall denote a constant which can
change at each occurrence.

\newsection{Linear Estimates in ${\mathbb R}_+\times{\mathbb R}^n$}

To handle the existence results for Klein-Gordon equations or
Dirichlet-wave equations, we shall need estimates which follow
immediately from estimates in H\"ormander \cite{H}.   To be more
precise, if
\begin{equation}\label{2.0}\{\Gamma_{t,x}\}=\{\partial_t, \partial_x,
\Omega_{jk}: \, 0\le j<k\le n\}
\end{equation}
 then we need the
following estimate to handle the case where $n=3$.

\begin{proposition}\label{prop2.3} Suppose that
$u\in C^\infty({\mathbb R}\times {\mathbb R}^3)$  satisfies
$u(t,x)=0$, $t\le 2B$, where $B$ is a fixed positive constant.
Suppose also that $(\square_{{\mathbb R}^{1+3}} +\MU^2)u(t,x)=0$
for $|x|>t-B$, Then there is a constant depending only on $B$  so
that when $\mu\ge1$
%\begin{multline}\label{2.3}
%\sup t^{n/2}|u(t,x)|\le C\sum_{|\alpha|\le
%5}\Bigl(\sum_{k=1}^\infty \sup_{\tau\in [2^{k-1},2^k]\cap
%[0,t]}2^k\|\Gamma^\alpha (\square+\MU^2)u(s,\cd)\|_2
%\\
%+ \sup_{\tau\in [0,1]}\|\Gamma^\alpha (\square
%+\MU^2)u(\tau,\cd)\|_2\Bigr).
%\end{multline}
\begin{equation}\label{2.3}
\sup_x t^{3/2}|u(t,x)|\le C\sum_{|\alpha|\le 5}\sum_{k}
\sup_{\tau\in [2^{k-1},2^{k+1}]\cap [2B,t]}2^k\|\Gamma^\alpha_{t,x}
(\square_{{\mathbb R}^{1+3}}+\MU^2)u(\tau,\cd)\|_2
\end{equation}
\end{proposition}

This follows exactly from the proof of Proposition 7.3.6 in
H\"ormander \cite{H} if one uses a variation of Lemma 7.3.4 there.
The relevant version for \eqref{2.3} is that if $w''+\MU^2w=h$ in
$[a,b]\subset {\mathbb R}$, then
$$\sup_{a\le \rho\le b}|w(\rho)|\le |w(a)|+|w'(a)|+\frac1{\mu}\int_a^b
|h(\rho)|d\rho.$$

In a similar manner one obtains the following analog of
Proposition 7.3.7 in \cite{H} for the case where $n\ge4$.

\begin{proposition}\label{prop2.4}Suppose that $n\ge4$ and that
$u\in C^\infty({\mathbb R}\times {\mathbb R}^n)$  satisfies
$u(t,x)=0$, $t\le 2B$, where $B$ is a fixed positive constant.
Suppose also that $(\square_{{\mathbb R}^{1+n}} +\MU^2)u(t,x)=0$
for $|x|>t-B$. Then if $n=4$ there is a constant depending only on
$B$  so that when $\mu\ge1$
%\begin{multline}\label{2.4}
%\sup t^2|u(t,x)| \le C\sum_{|\alpha|\le 7}\Bigl(\sum_k
%\sup_{\tau\in [2^{k-1},2^k]\cap [0,t]}(1+k)2^k\|\Gamma^\alpha
%(\square+\MU^2)u(s,\cd)\|_2
%\\
%+ \sup_{\tau\in [0,1]}\|\Gamma^\alpha (\square
%+\MU^2)u(\tau,\cd)\|_2\Bigr).
%\end{multline}
\begin{equation}\label{2.4}
\sup_x t^2|u(t,x)| \le C\sum_{|\alpha|\le 7}\sum_k \sup_{\tau\in
[2^{k-1},2^{k+1}]\cap [2B,t]}(1+|k|)2^k\|\Gamma_{t,x}^\alpha
(\square_{{\mathbb R}^{1+4}}+\MU^2)u(\tau,\cd)\|_2
\end{equation}
If $n\ge5$ there is a constant depending only on $B$ and $n$ so
that when $\mu\ge 1$
%\begin{multline}\label{2.5}
%\sup t^{1+\tfrac{n}4}|u(t,x)| \le C\sum_{|\alpha|\le
%n+3}\Bigl(\sum_k \sup_{\tau\in [2^{k-1},2^k]\cap
%[0,t]}2^k\|\Gamma^\alpha (\square+\MU^2)u(s,\cd)\|_2
%\\
%+ \sup_{\tau\in [0,1]}\|\Gamma^\alpha (\square
%+\MU^2)u(\tau,\cd)\|_2\Bigr).
\begin{equation}\label{2.5}
\sup_x t^{1+\tfrac{n}4}|u(t,x)| \le C\sum_{|\alpha|\le n+3}\sum_k
\sup_{\tau\in [2^{k-1},2^{k+1}]\cap [2B,t]}2^k\|\Gamma_{t,x}^\alpha
(\square_{{\mathbb R}^{1+n}}+\MU^2)u(\tau,\cd)\|_2
\end{equation}
\end{proposition}

To prove Theorem \ref{theorem1.3} (i.e., Neumann boundary
conditions with $m=0$), we shall also require weighted
$L^2_tL^2_x$ estimates of the type that were first used in
\cite{KSS2}.  Here and throughout, we will use the notation $\langle x\rangle = \langle r\rangle = \sqrt{1+|x|^2}$.

\begin{proposition}\label{prop2.5}  Fix $n\ge3$.  Then there is a uniform constant
$C$ which is independent of $\mu\ge 0$ so that if $u\in
C^\infty({\mathbb R}_+\times {\mathbb R}^n)$ vanishes for $t\le0$
and vanishes for large $x$ for every fixed $t$ then
\begin{multline}\label{2.6}
 \left(\int_0^T\int_{{\mathbb R}^n}\langle x\rangle^{-1}
\bigl[\, | \nabla_xu(s,x)\, |^2 + (1+ \mu)^{-1}\langle
x\rangle^{-1}\bigl(\, | \mu u(s,x)|^2 +|\partial_tu(s,x)|^2 \bigr)
\bigr]\, dx ds\right)^{1/2}
\\
\le C\bigl( \log (2+T)\bigr)^{1/2}\int_0^T\|(\square_{{\mathbb
R}^{1+n}} +\MU^2)u(s,\cd)\|_2\, ds.
\end{multline}
Additionally, if $\sigma>0$ is fixed there is a constant which is
independent of $\mu$ and $T$ so that
\begin{multline}\label{2.7}
 \left(\int_0^T\int_{{\mathbb R}^n}\langle x\rangle^{-1-\sigma}
\bigl[\, | \nabla_xu(s,x)\, |^2 +  (1+\mu)^{-1}\langle
x\rangle^{-1}\bigl(\, | \mu u(s,x)|^2 +|\partial_tu(s,x)|^2 \bigr)
\bigr]\, dx ds\right)^{1/2}
\\
\le C\int_0^T\|(\square_{{\mathbb R}^{1+n}} +\MU^2)u(s,\cd)\|_2\,
ds.
%\begin{multline}\label{2.7}
%\left(\int_0^T\int_{{\mathbb R}^n} \bigl(|\, \langle
%x\rangle^{-\sigma-1/2} u'(s,x)\, |^2 + |\, \langle
%x\rangle^{-\sigma-3/2} u(s,x)|^2 \bigr)\, dx ds\right)^{1/2}
%\\
%\le C\int_0^T\|(\square_{{\mathbb R}^{1+n}} +\MU^2)u(s,\cd)\|_2\,
%ds.
\end{multline}
\end{proposition}

To prove this we note that by making a dyadic decomposition of
$\{x\in {\mathbb R}^n: \, |x|\ge1\}$ and applying a scaling
argument, one sees that this result follows from the following
analog of Lemma 2.2 of \cite{SS}.

\begin{lemma}\label{lemma2.6}  Let $\beta(x)$ be a fixed smooth
function that is supported in $\{x\in {\mathbb R}^n: \, |x|\le
1\}$. Assume as above that $n\ge3$ is fixed.  Then there is a
constant $C$ which is independent of $\mu \ge0$ so that
\begin{multline}\label{2.8}
\int_{-\infty}^\infty\left( \bigl\|\nabla_x\beta(\cd)\bigl(
e^{it\sqrt{\MU^2-\Delta}}f\bigr)(t,\cd)\|^2_{L^2({\mathbb R}^n)}
+(1+\mu) \bigl\|\beta(\cd)\bigl(
e^{it\sqrt{\MU^2-\Delta}}f\bigr)(t,\cd)\|^2_{L^2({\mathbb
R}^n)}\right)\, dt
\\
\le C\int_{{\mathbb R}^n}|\hat f(\xi)|^2 \, (\MU^2+|\xi|^2)\,
d\xi.
\end{multline}
\end{lemma}

\noindent{\bf Proof of Lemma \ref{lemma2.6}:}   We just follow the
proof of the special case of $\gamma=1$ of Lemma 2.2 in \cite{SS}.

We first notice that by Plancherel's theorem over $t,x$ the left
side of \eqref{2.8} can be written as $(2\pi)^{n+1}$ times
$$\int_\mu^\infty \int\left| \int \hat \beta(\xi-\eta) \, \hat
f(\eta) \, \delta (\tau-\sqrt{\MU^2+|\eta|^2})\, d\eta\right|^2 \,
(1+\mu+|\xi|^2)\, d\xi d\tau.$$ We next apply the Schwarz
inequality over $\eta$ to bound this by
\begin{multline*}
\int_{\mu}^\infty \int (1+\mu+|\xi|^2) \left[ \int |\hat
\beta(\xi-\eta)|
 \delta(\tau-\sqrt{\MU^2+|\eta|^2}) \, d\eta\right]
 \\
\times
 \left[\int |\hat \beta(\xi-\eta)| \, |\hat f(\eta)|^2
\delta(\tau-\sqrt{\MU^2+|\eta|^2}) d\eta\right] \, d\xi d\tau.
\end{multline*}
Since
$\delta(\tau-\sqrt{\mu^2+r^2})=\frac{\tau}{\sqrt{\tau^2-\mu^2}} \,
\delta(r-\sqrt{\tau^2-\mu^2})$, using polar coordinates
$\eta=r\omega$, $r>0$, $\omega\in S^{n-1}$, we see that if
$\tau>\mu$ then
\begin{align*}
&(1+\mu+|\xi|^2)\int |\hat
\beta(\xi-\eta)|\delta(\tau-\sqrt{\mu^2+|\eta|^2}) \, d\eta
\\
&=c_n(1+\mu+|\xi|^2)\int_{S^{n-1}}|\hat
\beta(\xi-\sqrt{\tau^2-\mu^2}\, \omega)|
\tau(\sqrt{\tau^2-\mu^2})^{n-2}\, d\omega
\\
&\le C_N(1+\mu+|\xi|^2)(1+|\, |\xi|-\sqrt{\tau^2-\mu^2}\,
|)^{-N}\tau(\sqrt{\tau^2-\mu^2})^{n-2}(1+\sqrt{\tau^2-\mu^2}\,
)^{-(n-1)},
\end{align*}
for any $N$.

We conclude that there is a $C$ which is independent of $\mu$ so
that the left side of \eqref{2.8} is dominated by
\begin{multline*}
\iint |\hat\beta(\xi-\eta)|\hat f(\eta)|^2
(1+\mu+|\eta|^2)\sqrt{\mu^2+|\eta|^2}|\eta|^{n-2}
(1+|\eta|)^{-(n-1)}\, d\eta\,d\xi
\\
\le C\int |\hat f(\eta)|^2 (\mu^2+|\eta|^2)\,
d\eta,\end{multline*}
 using the assumption that
$n\ge3$ in the last step.
 \qed

\noindent{\bf Proof of Proposition \ref{prop2.5}:}
We will only prove \eqref{2.6}.  Similar techniques can be used to obtain \eqref{2.7}; see, e.g., \cite{Metcalfe}.

First note that for $|x|\ge T$, \eqref{2.6} follows from the
standard energy inequality. Next, we observe that when the norm in
the left side is taken over $\{|x|<1\}$, \eqref{2.6} follows from
\eqref{2.8}.  Indeed, by \eqref{2.8} and Duhamel's formula, we
have
\begin{equation}\label{2.9}
 \|\, |\nabla_x u| + (1+\mu)^{-1/2}(\mu|u|+|\partial_tu|)\, \|_{L^2([0,T]\times \{|x|<1\})} \le
C \int_0^T \|(\Box_{\R^{1+n}}+\mu^2)u(s,\cd)\|_2\:ds.\end{equation}

Let $u_R(t,x)=R^{-2}u(Rt,Rx)$, then
$\bigl((\square+\mu^2)u\bigr)(Rt,Rx)=(\square+(R\mu)^2)u_R(t,x)$.
So if we apply \eqref{2.9} with $u$ replaced by $u_R$ and $\mu$
replaced by $R\mu$ then we can change variables to see that the
analog of (2.5) holds without the log factor if in the left the
norms are taken over $[0,T]\times \{|x|\approx R\}$.  If we square
both sides and then sum  over dyadic $R=2^k<T$, we get (2.5). \qed

%Let $u_R(t,x)=u(Rt,Rx)$.  Notice that
%\begin{multline*}
%\|\langle x\rangle^{-1/2} u'\|_{L^2([0,T]\times \{x\in [R/2,R]\})} +
%\|\langle x\rangle^{-3/2} u\|_{L^2([0,T]\times \{x\in [R/2,R]\})}
%\\\le C R^{(n-2)/2} \Bigl(\|u'_R\|_{L^2([0,T/R]\times \{|x|<1\})} + \|u_R\|_{L^2([0,T/R]\times \{|x|<1\})}\Bigr)
%\end{multline*}
%Since $|(\Box_{\R^{1+n}}+\mu^2)u_R(t,x)|
%\le R^2 |(\Box_{\R^{1+n}}+\mu^2)u(Rt,Rx)|$ for $R>1$, we can apply \eqref{2.9} to see that this is
%$$\le C R^{(n+2)/2} \int_0^{T/R} \|(\Box_{\R^{1+n}}+\mu^2) u(Rt,Rx)\|_2\:ds
%= C \int_0^T \|(\Box_{\R^{1+n}}+\mu^2)u(s,\cd)\|_2\:ds.$$
%Thus, we have just shown that for $R>1$
%\begin{multline}\label{2.10}
%\|\langle x\rangle^{-1/2} u'\|_{L^2([0,T]\times \{x\in [R/2,R]\})} + \|\langle x\rangle^{-3/2} u\|_{L^2([0,T]
%\times \{x\in [R/2,R]\})} \\\le C \int_0^T \|(\Box_{\R^{1+n}}+\mu^2) u(s,\cd)\|_2\:ds.
%\end{multline}
%
%In order to complete the proof, we need to show that the bound holds when the spatial norm on the left
%side of \eqref{2.6} is taken over $1<|x|<t$.  By squaring the left side of \eqref{2.6} and then summing over
%the dyadic pieces with $1<R=2^j<t$, \eqref{2.10} gives \eqref{2.6}.
%\qed

\noindent{\bf Remark:}  Although we shall not use it, the above
proof yields a somewhat weaker result when $n=2$.  Namely, in this
case one has that the analogs of \eqref{2.6} and \eqref{2.7} hold
with a constant that is independent of $\mu\ge1$ (instead of
$\mu\ge0$ as in the $n\ge3$ case).

\medskip

To handle Neumann boundary conditions when $m=0$ in \eqref{1.1},
we also require the following well known Klainerman-Sobolev
estimates.  See \cite{K}.

\begin{proposition}\label{prop2.2}  Let $n\ge2$, and set
$$\{Z \}=\{\partial_x, \Omega_{jk}: \, 1\le j<k\le n\}.$$
Then if $R\ge1$
\begin{equation}\label{2.2}
R^{(n-1)/2}\sup_{R-1\le |x|\le R}|h(t,x)|\le C\sum_{|\alpha|\le
(n+2)/2}\|Z^\alpha h(t,\cd)\|_{L^2(R-2\le |x|\le R+1)}.
\end{equation}
\end{proposition}

%%%%%%%%%%%%%%%%%%%%%%%%%%%%%%%%%%%%%%%%%%%%%%%%%%%%%%%%
\newsection{Linear estimates for  waveguides ${\mathbb
R}_+\times ({\mathbb R}^n\times \Omega)$}

The purpose of this section is to show that the estimates for
solutions of wave and Klein-Gordon equations from the previous
section carry over to the setting of waveguides if one gives up a
few derivatives, depending on the dimension $d$ of $\Omega \subset
{\mathbb R}^d$.

To prove these results we need to recall a few basic facts from
the spectral theory of $\Delta_\Omega$ and from its elliptic
regularity theory (see, e.g., \cite{Taylor} for proofs;  see also \cite{gilbarg}).  As
before, $\Delta_\Omega$ denotes either the Dirichlet Laplacian, in
which case we impose boundary conditions
\begin{equation}\label{3.1}
h|_{\partial \Omega}=0,
\end{equation}
or the Neumann Laplacian where we require
\begin{equation}\label{3.2}
\partial_\nu h|_{\partial \Omega}=0
\end{equation}
with $\partial_\nu$ denoting the normal derivative.

We first recall that since we are assuming that $\Omega$ is
compact with smooth boundary, the spectrum of $-\Delta_\Omega$ is
discrete and nonnegative.  We let $\lambda_1^2\le \lambda_2^2\le
\lambda_3^2\le \cdots$ denote the eigenvalues counted with respect
to multiplicity.  Then
\begin{equation}\label{3.3}
0<\lambda_1^2 \le \lambda^2_2\le \cdots, \quad \text{for the
Dirichlet Laplacian},
\end{equation}
and
\begin{equation}\label{3.4}
0=\lambda_1^2<\lambda_2^2\le \lambda_3^2\cdots, \quad \text{for
the Neumann Laplacian}.
\end{equation}
For either of these Laplacians, we shall let $E_j: \,
L^2(\Omega)\to L^2(\Omega)$ denote the projection onto the $j$th
eigenspace.  Thus, if $h\in L^2(\Omega)$, we have that $E_jh$ is
smooth,
\begin{equation}\label{3.5}
-\Delta_\Omega E_jh(x)=\lambda_j^2 E_jh(x),
\end{equation}
and we have Plancherel's theorem
\begin{equation}\label{3.6}
\|h\|_{L^2(\Omega)}^2 = \sum_{j=1}^\infty \|E_jh\|_2^2.
\end{equation}
Also, by the well known Weyl formula we have $\lambda_j\approx
j^{1/d}$, $j=2,3,\dots$.  Therefore, if $h\in
C^\infty(\overline{\Omega})$,
\begin{multline}\label{3.77}
(1+j)^{2/d}\|E_jh\|_{L^2(\Omega)}\le
C\|(I-\Delta_\Omega)E_jh\|_{L^2(\Omega)}
\\
= C\|E_j(I-\Delta_\Omega)h\|_{L^2(\Omega)}, \quad j=1,2,3,\dots,
\end{multline}
assuming that either \eqref{3.1} or \eqref{3.2} holds.

Having reviewed the spectral theory, let us recall the basic
elliptic regularity estimate:

\begin{lemma}\label{lemma3.1}  Suppose that $h\in
C^\infty(\overline{\Omega})$ and that we have the boundary
conditions \eqref{3.1} or \eqref{3.2}.  Then for $N=2,3,\dots$
there is a constant $C=C_{N,\Omega}$ so that
$$\sum_{|\alpha|\le N}\|\partial^\alpha_y h\|_{L^2(\Omega)}
\le C\sum_{|\alpha|\le N-2}\|\partial_y^\alpha \Delta_\Omega
h\|_{L^2(\Omega)}+C\|h\|_{L^2(\Omega)}.$$
\end{lemma}

Using this result we can obtain the following

\begin{lemma}\label{lemma3.2}  Suppose that $u(t,x,y)\in
C^\infty({\mathbb R}_+\times \Rn \times \Omega)$ and that either
$u(t,x,y)=0$ for all $y\in \partial \Omega$ or $\partial_\nu
u(t,x,y)=0$ for all $y\in \partial\Omega$.  Then if $m\ge 0$ and
$N=2,4,6,\dots$ are fixed
\begin{multline*}
\sum_{|\alpha|\le N}\|\partial^\alpha_y
u(t,x,\cd)\|_{L^2(\Omega)}
\\
\le C\sum_{|\alpha| \le
N}\|\partial^\alpha_{t,x}u(t,x,\cd)\|_{L^2(\Omega)}+C\sum_{|\alpha|\le
N-2}\|\partial^\alpha_{t,x,y}(\square+m^2)u(t,x,\cd)\|_{L^2(\Omega)}.
\end{multline*}
Moreover, for any $N=2,3,4,\dots,$
\begin{multline*}
\sum_{|\alpha|\le N} \|\partial^\alpha_y u(t,x,\cd)\|_{L^2(\Omega)}
\\\le C \sum_{\substack{|\alpha|+|\beta|\le N\\|\beta|\le 1}}
\|\partial^\alpha_{t,x}\partial_y^\beta u(t,x,\cd)\|_{L^2(\Omega)} +
C\sum_{|\alpha|\le N-2} \|\partial_{t,x,y}^\alpha (\square +m^2)u(t,x,\cd)\|_{L^2(\Omega)}.
\end{multline*}
\end{lemma}

The proof is quite simple.  Since $-\Delta_\Omega = (\square +
m^2)-\partial^2_t+\Delta_{\Rn}-m^2$, and since $\partial_t$ and
$\nabla_x$ preserve the boundary conditions, one just uses Lemma
\ref{lemma3.1} and an induction argument.

We also require a simple lemma which shows how the spectrum of
$\Delta_\Omega$ is relevant for our equations involving
$(\square+m^2)$.

\begin{proposition}\label{prop3.3}  Let $\Delta_\Omega$ denote either the
Dirichlet Laplacian or the Neumann Laplacian on $\Omega$.  Then,
if $u\in C^\infty({\mathbb R}_+\times \Rn\times \Omega)$ satisfies
the relevant boundary condition \eqref{1.3} or \eqref{1.4}, it
follows that
$$E_j (\square+m^2) u(t,x,y)=(\partial_t^2-\Delta_{\Rn}+m^2+\lambda_j^2)E_ju(t,x,y),$$
and so
$$\sum_j\|(\partial_t^2-\Delta_{\Rn}+m^2+\lambda_j^2)E_ju(t,x,\cd)\|_{L^2(\Omega)}^2=
\|(\square+m^2)u(t,x,\cd)\|_{L^2(\Omega)}^2.
$$
\end{proposition}

This follows immediately from \eqref{3.5} and \eqref{3.6}.

Next we need the following consequence of Sobolev's lemma for
$\Omega$ and the above elliptic regularity estimates.

\begin{proposition}\label{prop3.4}  Let $u$ be as above.  Then if
$m\ge0$ is fixed
\begin{multline}\label{3.7}
|\partial^\beta_y u(t,x,y)|\le C\sum_{|\alpha|\le
|\beta|+(d+4)/2}\|\partial^\alpha_{t,x} u(t,x,\cd)\|_{L^2(\Omega)}
\\
+C\sum_{|\alpha|\le
|\beta|+d/2}\|\partial^\alpha_{t,x,y}(\square+m^2)u(t,x,\cd)\|_{L^2(\Omega)}.
\end{multline}
\end{proposition}

The proof is simple.  One first uses Sobolev's lemma for $\Omega$
to obtain
$$|\partial^\beta_y u(t,x,y)|\le C\sum_{|\alpha|\le
|\beta|+(d+2)/2}\|\partial^\alpha_y u(t,x,\cd)\|_{L^2(\Omega)}.$$
Therefore, if we apply Lemma \ref{lemma3.2} (and sum over indices of
length one more if $|\beta|+(d+2)/2$ is odd), we obtain
\eqref{3.7}.

We can now state one of our most important estimates.  We shall
let
\begin{equation}\label{3.8}
\{\Gamma\}= \{\Gamma_{t,x}\}\cup \{\partial_y\}=
\{\partial_t,\partial_x, \Omega_{jk}, \partial_y: \, 0\le j< k \le
n\}.
\end{equation}

%\begin{proposition}\label{prop3.5}  Fix $B$ and suppose that $u\in
%C^\infty({\mathbb R}_+\times \Rn\times \overline{\Omega})$
%satisfies $u(t,x,y)=0$, $t\le 2B$ and $u(t,x,y)=0$ if $|x|>t-B$.
%Suppose also that either $m\ge 0$ and $u(t,x,y)=0$, $y\in \partial
%\Omega$ or that $m>0$ and $\partial_\nu u(t,x,y)=0$, $y\in
%\partial \Omega$.  Then, if $n=3$,
%\begin{multline*}
%t^{3/2} |\Gamma^\beta u(t,x,y)| \le C\sum_{|\alpha|\le |\beta| +
%5+(5d+3)/2}\sum_k \sup_{\tau\in [2^{k-1},2^k]\cap[2B,t]}2^k
%\|\Gamma^\alpha (\square+m^2)u(s,\cd)\|_2
%\\ +Ct^{3/2}\sum_{|\alpha|\le
%|\beta|+3+(d-1)/2}\|\Gamma^\alpha(\square+m^2)u(t,\cd)\|_2.
%\end{multline*}
%If $n=4$,
%\begin{multline*}t^2 |\Gamma^\beta u(t,x,y)|\le
%C\sum_{|\alpha|\le |\beta|+7+(5d+3)/2}\sum_k \sup_{\tau\in
%[2^{k-1},2^k]\cap [2B,t]} (1+k)2^k\|\Gamma^\alpha
%(\square+m^2)u(s,\cd)\|_2
%\\
%+Ct^2\sum_{|\alpha|\le
%|\beta|+5+(d-1)/2}\|\Gamma^\alpha(\square+m^2)u(t,\cd)\|_2.
%\end{multline*}
% Finally, if $n\ge5$
%\begin{multline*}t^{1+\tfrac{n}4} |\Gamma^\beta u(t,x,y)|\le
%C\sum_{|\alpha|\le |\beta|+n+3+(5d+3)/2}\sum_k \sup_{\tau\in
%[2^{k-1},2^k] \cap [2B,t]}2^k \|\Gamma^\alpha (\square
%+m^2)u(s,\cd)\|_2
%\\
%+Ct^{1+\tfrac{n}4}\sum_{|\alpha|\le
%|\beta|+n+1+(d-1)/2}\|\Gamma^\alpha(\square+m^2)u(t,\cd)\|_2.
%\end{multline*}
% Here as above, we are assuming that $d$
%is the dimension of $\Omega$.
%\end{proposition}

\begin{proposition}\label{prop3.5}  Fix $B$ and suppose that $u\in
C^\infty({\mathbb R}_+\times \Rn\times \overline{\Omega})$
satisfies $u(t,x,y)=0$, $t\le 2B$ and $(\Box+m^2)u(t,x,y)=0$ if $|x|>t-B$.
Suppose also that either $m\ge 0$ and $u(t,x,y)=0$, $y\in \partial
\Omega$ or that $m>0$ and $\partial_\nu u(t,x,y)=0$, $y\in
\partial \Omega$.  Then, if $n=3$,
\begin{multline*}
t^{3/2} |\Gamma^\beta u(t,x,y)| \le C\sum_{|\alpha|\le |\beta| +
5+(5d+4)/2}\sum_k \sup_{\tau\in [2^{k-1},2^{k+1}]\cap[2B,t]}2^k
\|\Gamma^\alpha (\square+m^2)u(\tau,\cd)\|_2
\\ +Ct^{3/2}\sum_{|\alpha|\le
|\beta|+(5d+5)/2}\|\Gamma^\alpha(\square+m^2)u(t,\cd)\|_2.
\end{multline*}
If $n=4$,
\begin{multline*}t^2 |\Gamma^\beta u(t,x,y)|\\\le
C\sum_{|\alpha|\le |\beta|+7+(5d+4)/2}\sum_k \sup_{\tau\in
[2^{k-1},2^{k+1}]\cap [2B,t]} (1+|k|)2^k\|\Gamma^\alpha
(\square+m^2)u(\tau,\cd)\|_2
\\
+Ct^2\sum_{|\alpha|\le
|\beta|+(5d+n+2)/2}\|\Gamma^\alpha(\square+m^2)u(t,\cd)\|_2.
\end{multline*}
 Finally, if $n\ge5$
\begin{multline*}t^{1+\tfrac{n}4} |\Gamma^\beta u(t,x,y)|\\\le
C\sum_{|\alpha|\le |\beta|+n+3+(5d+4)/2}\sum_k \sup_{\tau\in
[2^{k-1},2^{k+1}] \cap [2B,t]}2^k \|\Gamma^\alpha (\square
+m^2)u(\tau,\cd)\|_2
\\
+Ct^{1+\tfrac{n}4}\sum_{|\alpha|\le
|\beta|+(5d+n+2)/2}\|\Gamma^\alpha(\square+m^2)u(t,\cd)\|_2.
\end{multline*}
 Here as above, we are assuming that $d$
is the dimension of $\Omega$.
\end{proposition}

%This follows immediately from Propositions
%\ref{prop2.3}-\ref{prop2.5} and Proposition \ref{prop3.4}. One
%uses the fact that $\{\partial_t,\partial_x,\Omega_{jk}, \, 0\le
%j<k\le n\}$ preserve the boundary conditions.

\noindent{\bf Proof of Proposition \ref{prop3.5}:}
To prove this, we first note that since $\{\Gamma_{t,x}\}$ commute
with $\square +m^2$ and preserve the boundary conditions, we only
have to prove the above estimates when
$\Gamma^\beta=\partial_y^\beta$, where, as before $y$ is the
$\Omega$-variable.  If we apply Proposition \ref{prop3.4} and
orthogonality, we conclude that
\begin{align*}
|\partial^\beta_yu(t,x,y)|^2&\le C\sum_{|\alpha|\le
|\beta|+(d+4)/2}
\|\partial^\alpha_{t,x}u(t,x,\cd)\|^2_{L^2(\Omega)}
\\
&\qquad\qquad\qquad
 + C\sum_{|\alpha|\le
|\beta|+d/2}\|\partial^\alpha_{t,x,y}(\square+m^2)u(t,x,\cd)\|_{L^2(\Omega)}^2
\\
&\le C\sum_{|\alpha|\le |\beta|+(d+4)/2}\sum_{j=1}^\infty
\|E_j\partial^\alpha_{t,x}u(t,x,\cd)\|^2_{L^2(\Omega)}
\\
&\qquad\qquad\qquad
 + C\sum_{|\alpha|\le
|\beta|+(d+n+2)/2}\|\partial^\alpha_{t,x,y}(\square+m^2)u(t,\cd)\|_{L^2(\Rn\times\Omega)}^2,
\end{align*}
using Sobolev's lemma for $\Rn$ in the last step.

To estimate the first term in the right, we use \eqref{3.77} to
get
\begin{align*}
(1+j)^{2/d}\sum_{|\alpha|\le
|\beta|+(d+4)/2}\|E_j&\partial^\alpha_{t,x}u(t,x,\cd)\|_{L^2(\Omega)}
\\
&\le C\sum_{|\alpha|\le
|\beta|+(d+4)/2}\|E_j(I-\Delta_\Omega)\partial^\alpha_{t,x}u(t,x,\cd)\|_{L^2(\Omega)}
\\
&\le C\sum_{|\alpha|\le 2+
|\beta|+(d+4)/2}\|E_j\partial^\alpha_{t,x}u(t,x,\cd)\|_{L^2(\Omega)}
\\
&\quad+ C\sum_{|\alpha|\le
|\beta|+(d+4)/2}\|E_j\partial_{t,x}^\alpha(\square+m^2)u(t,x,\cd)\|_{L^2(\Omega)}.
\end{align*}
By induction,
\begin{align*}
\sum_{|\alpha|\le
|\beta|+(d+4)/2}\|E_j\partial^\alpha_{t,x}u(t,x,\cd)\|_{L^2(\Omega)}&\le
\frac{C}{(1+j)^2}\sum_{|\alpha|\le
|\beta|+(5d+4)/2}\|E_j\partial^\alpha_{t,x}u(t,x,\cd)\|_{L^2(\Omega)}
\\
&\quad +C\sum_{|\alpha|\le
|\beta|+5d/2}\|E_j\partial^\alpha_{t,x}(\square+m^2)u(t,x,\cd)\|_{L^2(\Omega)}.
\end{align*}
Therefore, by \eqref{3.6} and Sobolev's lemma for $\Rn$ we
conclude that
\begin{multline*}
\sum_{|\alpha|\le
|\beta|+(d+4)/2}\sum_{j=1}^\infty\|E_j\partial^\alpha_{t,x}u(t,x,\cd)\|_{L^2(\Omega)}^2
\\\le C\sum_{|\alpha|\le
|\beta|+(5d+n+2)/2}\|\partial^\alpha_{t,x}(\square+m^2)u(t,\cd)\|^2_{L^2(\Rn\times
\Omega)}
\\
+C\sum_{j=1}^\infty\frac{1}{(1+j)^4}\Bigl(\sum_{|\alpha|\le
|\beta|+(5d+4)/2}\|E_j\partial^\alpha_{t,x}u(t,x,\cd)\|_{L^2(\Omega)}^2\Bigr).
\end{multline*}
Thus, if we combine these two steps, we conclude that
\begin{align}\label{3.100}
t^n|\partial^\beta_yu(t,x,y)|^2&\le Ct^n\sum_{|\alpha|\le
|\beta|+(5d+n+2)/2}\|\partial^\alpha_{t,x}(\square+m^2)u(t,\cd)\|_{L^2(\Rn\times\Omega)}^2
\\
&\quad +Ct^n\sum_{j=1}^\infty
\frac1{(1+j)^4}\Bigl(\sum_{|\alpha|\le
|\beta|+(5d+4)/2}\|E_j\partial^\alpha_{t,x}u(t,x,\cd)\|_{L^2(\Omega)}^2\Bigr).
\notag
\end{align}

We can use Proposition \ref{prop3.3} and Propositions
\ref{prop2.3} and \ref{prop2.4} to estimate each of the summands
in the last term since our assumptions regarding $m$ and the
boundary conditions ensure that $m^2+\lambda_j^2\ge c$, for some
constant $c>0$.  When $n=3$, we get from Proposition \ref{prop2.3}
that there is a constant $C$ independent of $j$ such that
\begin{align*}
t&^{3/2}\sum_{|\alpha|\le
|\beta|+(5d+4)/2}\|E_j\partial^\alpha_{t,x}u(t,x,\cd)\|_{L^2(\Omega)}
\\
&\le C\sum_{|\alpha|\le |\beta|+5+(5d+4)/2}\sum_k\sup_{\tau\in
[2^{k-1},2^{k+1}]\cap[2B,t]}2^k\\&\quad\quad\quad\quad\quad\quad
\quad\quad\quad\quad\quad\quad\quad
 \times
\|\Gamma^\alpha_{t,x}(\partial_t^2-\Delta_{{\mathbb
R}^3}+m^2+\lambda_j^2)E_ju(\tau,\cd)\|_{L^2({\mathbb
R}^3\times\Omega)}
\\
&\le C\sum_{|\alpha|\le |\beta|+5+(5d+4)/2}\sum_k\sup_{\tau\in
[2^{k-1},2^{k+1}]\cap[2B,t]}2^k\|\Gamma^\alpha_{t,x}(\square+m^2)u(\tau,\cd)\|_{L^2({\mathbb
R}^3\times\Omega)},
\end{align*}
using Proposition \ref{prop3.3} in the last step.  This inequality
and \eqref{3.100} yield the estimate in Proposition \ref{prop3.5}
for $n=3$.  The estimates for $n\ge4$ follow from a similar
argument.  \qed

%Note that our assumptions regarding $m$ and the boundary
%conditions ensure that $m^2+\lambda_j^2\ge c$, for some constant
%$c>0$.  Using this fact, we can use Proposition \ref{prop3.3} to
%estimate the first term in the right side of the last inequality.
%If $n=3$ we get
%\begin{align*}
%t^{3/2}&\sum_{|\alpha|\le
%|\beta|+(d+3)/2}\|E_j\partial^\alpha_{t,x}u(t,x,\cd)\|_{L^2(\Omega)}
%\\
%&\le C\sum_{|\alpha|\le |\beta|+5+(d+3)/2}\sum_k\sup_{\tau \in
%[2^{k-1},2^k]\cap[2B,t]}2^k\|\Gamma^\alpha_{t,x}(\partial^2_t-\Delta_{{\mathbb
%R}^3}+m^2+\lambda_j^2)E_ju(t,\cd)\|_{L^2({\mathbb R}^3\times
%\Omega)}
%\\
%&\le C\sum_{|\alpha|\le |\beta|+5+(d+3)/2}\sum_k\sup_{\tau \in
%[2^{k-1},2^k]\cap[2B,t]}2^k\|E_j\Gamma^\alpha_{t,x}(\square+m^2)u(t,\cd)\|_{L^2({\mathbb
%R}^3\times \Omega)}
%\\
%&\le \frac{C}{(1+j)^2}\sum_{|\alpha|\le
%|\beta|+5+(d+3)/2}\sum_k\sup_{\tau \in
%[2^{k-1},2^k]\cap[2B,t]}2^k\|(I-\Delta_\Omega)^d
%E_j\Gamma^\alpha_{t,x}(\square+m^2)u(t,\cd)\|_{L^2({\mathbb
%R}^3\times \Omega)}
%\\
%&\le \frac{C}{(1+j)^2}\sum_{|\alpha|\le
%|\beta|+5+(5d+3)/2}\sum_k\sup_{\tau \in
%[2^{k-1},2^k]\cap[2B,t]}2^k\|\Gamma^\alpha(\square+m^2)u(t,\cd)\|_{L^2({\mathbb
%R}^3\times \Omega)},
%\end{align*}
%using \eqref{3.77} in the last two steps.  If we combine the last
%two inequalities we get the estimate in Proposition \ref{prop3.5}
%for $n=3$.  The estimates for $n\ge4$ follow from a similar
%argument.

To prove Theorem \ref{theorem1.3} we need the following

\begin{proposition}\label{prop3.6}  Suppose that $u\in
C^\infty({\mathbb R}_+\times \Rn \times \overline{\Omega})$
satisifes $u(t,x,y)=0$, $t\le 0$ and $\partial_\nu u(t,x,y)=0$,
$y\in \partial \Omega$.  Then if $n\ge3$
\begin{align}\label{3.10}
&\sum_{|\alpha|\le
N}\left(\|\Gamma^\alpha\nabla_{t,x}u(T,\cd)\|_{L^2(\Rn\times\Omega)}
+ (\log (2+T))^{-1/2}\|\, \langle x\rangle^{-1/2}\Gamma^\alpha
\nabla_xu\|_{L^2([0,T]\times \Rn\times\Omega)}\right)
\\
&\le C\int_0^T\sum_{|\alpha|\le N}\|\Gamma^\alpha \square
u(t,\cd)\|_{L^2(\Rn\times\Omega)}\, dt \notag
\\
%&\qquad
&+C\sum_{|\alpha|\le N-1}\Bigl(\|\Gamma^\alpha \square
u(T,\cd)\|_{L^2(\Rn\times \Omega)} \notag
%\\
%&\qquad\qquad\qquad
+(\log(2+T))^{-1/2}\|\, \langle x\rangle^{-1/2}\Gamma^\alpha
\square u\|_{L^2([0,T]\times \Rn\times\Omega)}\Bigr). \notag
\end{align}
\end{proposition}

\noindent{\bf Proof of Proposition \ref{prop3.6}:}
%To prove this,
Let us start out by showing that the second term in
the left side satisfies the desired estimates.  By Lemma
\ref{lemma3.2}, it suffices to prove the variant where all of the
$\Gamma^\alpha$ are of the form $\Gamma^\alpha_{t,x}$.  Since
these vector fields preserve the boundary conditions, we see that
it suffices to prove the estimate for $N=0$.  But this follows
from Proposition \ref{prop2.5} and Proposition \ref{prop3.3}:
\begin{align*}
(\log(2+T))^{-1}&\|\, \langle
x\rangle^{-1/2}\nabla_xu\|^2_{L^2([0,T]\times \Rn\times\Omega)}
\\
&=(\log(2+T))^{-1}\sum_{j=1}^\infty\|\, \langle
x\rangle^{-1/2}\nabla_xE_ju\|^2_{L^2([0,T]\times \Rn\times
\Omega)}
\\
&\le C\sum_{j=1}^\infty\left[
\int_0^T\|(\partial_t^2-\Delta_{\Rn}+\lambda_j^2)E_ju(t,\cd)\|_{L^2(\Rn\times\Omega)}\,
dt \right]^2
\\
&\le C\left[
\int_0^T\Bigl(\sum_{j=1}^\infty\|(\partial_t^2-\Delta_{\Rn}+\lambda_j^2)
E_ju(t,\cd)\|_{L^2(\Rn\times\Omega)}^2\Bigr)^{1/2}\, dt \right]^2
\\
&\le C\left[\int_0^T\| \square u(t,\cd)\|_{L^2(\Rn\times\Omega)}\,
dt\right]^2.
\end{align*}
This completes the proof that the second term in the left side of
\eqref{3.10} satisfies the desired bounds.  If one uses the energy
inequality for $\Rn$, then similar arguments imply that the other
term is also dominated by the right side.\qed

One can similarly use \eqref{2.7} to show that
\begin{multline}\label{3.12}
\sum_{|\alpha|\le N} \Bigl(\|\Gamma^\alpha \nabla_{t,x}
u(T,\cd)\|_{L^2(\R^n\times\Omega)} + \|\langle
x\rangle^{-1/2-\sigma} \Gamma^\alpha
\nabla_xu\|_{L^2([0,T]\times\R^n\times\Omega)}\Bigr)
\\\le C\int_0^T \sum_{|\alpha|\le N} \|\Gamma^\alpha \Box u(t,\cd)\|_{L^2(\R^n\times\Omega)}\:dt
\\+C\sum_{|\alpha|\le N-1}\Bigl(\|\Gamma^\alpha \Box u(T,\cd)\|_{L^2(\R^n\times\Omega)} + \|\langle x\rangle^{-1/2-\sigma}
\Gamma^\alpha \Box u\|_{L^2([0,T]\times\R^n\times\Omega)}\Bigr)
\end{multline}
for any $\sigma>0$.

%%%%%%%%%%%%%%%%%%%%%%%%%%%%%%%%%%%%%%%%%%%%%%%%%%%%%%%%%%%%%%%%%
\newsection{Existence theorems for Dirichlet boundary conditions}

Proving the existence results for \eqref{1.1} with Dirichlet
boundary conditions \eqref{1.3} is straightforward.  To apply
Proposition \ref{prop3.5}, as in \cite{H}, \cite{kkg}, one shifts
the time variable by $2B$ so that the initial condition is at
$t=2B$.  Then by well known local existence theorems (see e.g.
\cite{KSS}) if $\varepsilon>0$ is sufficiently small  and $N$ is large enough in
\eqref{1.6} there is always a solution of \eqref{1.1}, \eqref{1.3}
on $[2B,2B+1]$. To reduce to the case of zero initial data, we fix
$\eta\in C^\infty(\mathbb{R})$ so that $\eta(t)=1$, $t\le
2B+\tfrac12$ and $\eta(t)=0$, $t\ge 2B+1$.  Then
$$
u_0(t,x)=
\begin{cases}
\eta(t)u(t,x), \quad t\le 2B+1
\\
0, \quad \text{otherwise}
\end{cases}
$$
satisfies
$$(\square+m^2)u_0=\eta Q(u,u',u'')+[\square,\eta]u$$
with Dirichlet boundary conditions.  Therefore,
$(\square+m^2)u=Q(u,u',u'')$ for $2B<t<T$ and $(x,y)\in
\Rn\times\Omega$ if and only if $w=u-u_0$ solves
\begin{equation}\label{4.1}
\begin{cases}
(\square+m^2)w=(1-\eta)Q(u_0+w,(u_0+w)',(u_0+w)'')-[\square,\eta](u_0+w)
\\
w(t,x,y)=0, \quad y\in \partial\Omega
\\
w(t,x,y)=0, \quad t\le 2B.
\end{cases}\end{equation}
Since the local existence results imply that
%$$
\begin{equation}\label{4.2}
\sup_t \sum_{|\alpha|\le N}\|\Gamma^\alpha
u_0(t,\cd)\|_{L^2(\Rn\times \Omega)}\le C\varepsilon,
\end{equation}
%%$$
one can argue as in H\"ormander \cite{H} \S 7.4 that if $N$ is
large enough and if $\varepsilon$ is small enough then there is a
global solution of \eqref{4.1} if $n\ge3$.  One uses Proposition
\ref{prop3.5} and energy estimates.  The latter say that if
$\gamma^{jk}\in C^\infty$,
$$\sum_{j,k=0}^{n+d}|\gamma^{jk}|\le 1/2,$$
and if
$$\begin{cases}
(\square+m^2)w+\sum_{j,k=0}^{n+d}\gamma^{jk}(t,x,y)\partial_j\partial_kw
=F, \quad 2B\le t\le T
\\
w(t,x,y)=0, \quad y\in \partial\Omega
\\
w(t,x,y)=0, \quad t\le 2B,
\end{cases}
$$
then if $w$ vanishes for large $|x|$ we have
\begin{multline}\label{4.3}
\|\nabla_{t,x,y}w(t,\cd)\|_2+m\|w(t,\cd)\|_2
\\
\le 2\exp \Bigl(\int_0^t
2\sum_{i,j,k=0}^{n+d}\|\partial_i\gamma^{jk}(s,\cd)\|_\infty\,
ds\Bigr) \int_0^t\|F(s,\cd)\|_2\, ds.
\end{multline}
Here, as before, $x_0=t$, while $x_{n+j}=y_j$, $1\le j\le d$.
Since $w$ has Dirichlet boundary conditions the proof of Lemma
7.4.1 in \cite{H} shows that \eqref{4.3} holds.

In order to use the well-known arguments of Klainerman \cite{kkg} to
show global existence, we will require a version of \eqref{4.3} where
$w$ is replaced by $\Gamma^\alpha w$.  We begin with the case
$\Gamma=\partial_y$.  Here, since $\nabla_{t,x}$ preserves the
boundary condition, it follows from Lemma \ref{lemma3.2} that if
$|\alpha|=N$,
\begin{multline}\label{4.4}
\|\partial^\alpha_y \nabla_{t,x,y} w(t,\cd)\|_2 + m\|\partial^\alpha_y
w(t,\cd)\|_2
\\\le C \sum_{\substack{|\beta|+|\gamma|\le N+1\\|\gamma|\le 1}}
\|\partial_{t,x}^\beta \partial_y^\gamma w(t,\cd)\|_2
+C
\sum_{|\beta|\le N-1} \|\partial^\beta_{t,x,y} (\square+m^2)w(t,\cd)\|_2.
\end{multline}
Thus, since the vectors fields $\Gamma_{t,x}$ preserve the boundary
conditions and
\begin{multline*}
(\square+m^2)\Gamma^\alpha_{t,x} w + \sum_{j,k=0}^{n+d} \gamma^{jk}
\partial_j\partial_k\Gamma^\alpha_{t,x} w
\\=\Gamma^\alpha_{t,x} F +
\sum_{j,k=0}^{n+d}[\gamma^{jk},\Gamma_{t,x}^\alpha]\partial_j\partial_k
w + \sum_{j,k=0}^{n+d} \gamma^{jk}[\partial_j\partial_k,\Gamma_{t,x}^\alpha]w,
\end{multline*}
it follows easily from the proof of Lemma 7.4.1 in \cite{H}
 and \eqref{4.4} that for $N=0,1,2,\dots$,
\begin{multline}\label{4.5}
\sum_{|\alpha|\le N} \|\Gamma^\alpha \nabla_{t,x,y} w(t,\cd)\|_2 +
m\sum_{|\alpha|\le N} \|\Gamma^\alpha w(t,\cd)\|_2
\\\le C \exp \Bigl(\int_0^t
2\sum_{i,j,k=0}^{n+d} \|\partial_i\gamma^{jk}(s,\cd)\|_\infty\,
ds\Bigr) \sum_{|\alpha|\le N} \int_0^t \|\Gamma^\alpha_{t,x}
F(s,\cd)\|_2\, ds
\\+  C \exp \Bigl(\int_0^t
2\sum_{i,j,k=0}^{n+d} \|\partial_i\gamma^{jk}(s,\cd)\|_\infty
\,ds\Bigr)\sum_{|\alpha|\le N} \sum_{j,k=0}^{n+d} \int_0^t
\|[\gamma^{jk},\Gamma_{t,x}^\alpha]
\partial_j\partial_k w(s,\cd)\|_2\, ds
\\+C
\sum_{|\beta|\le N-1} \|\partial^\beta_{t,x,y} (\square+m^2)w(t,\cd)\|_2.
\end{multline}

Standard
arguments show that \eqref{4.5} and Proposition \ref{prop3.5}
imply the existence results for the nonlinear Dirichlet-wave
equations \eqref{1.1} when $m>0$.  Since Proposition \ref{prop3.5}
remains valid in the Dirichlet case when $m=0$, one also gets the
existence results for $m=0$.   This requires that
$\|w(t,\cd)\|_2$ is dominated by the right side of
\eqref{4.3} in the Dirichlet case
due to the fact that Poincar\'e's lemma gives
$\|w(t,\cd)\|_2\le C\|\nabla_y w(t,\cd)\|_2$.

%%%%%%%%%%%%%%%%%%%%%%%%%%%%%%%%%%%%%%%%%%%%%%%%%%%%%%%
\newsection{Existence theorems for Neumann boundary conditions}

%$$\square u=Q(u',u'')$$
%
%%\newsection{Global existence for nonlinear wave equations involving Dirichlet boundary
%%conditions when $n=3$}
%
%\newsection{Global existence for nonlinear wave equations  when $n\ge4$}
%
%$$\square u=Q(u',u'')$$
%
%Here get $u'=O(\varepsilon/t)$.  Use hyperbolic rotations and
%elliptic regularity when $|x|<t/2$ and use regular rotations for
%$|x|>t/2$

The proof of Theorem \ref{theorem1.2} follows from the same
arguments that were just used for the proof of Theorem
\ref{theorem1.1}.  Since one has Neumann boundary conditions,
though, one has to make an extra assumption in order for
\eqref{4.3} to hold for solutions of
\begin{equation}\label{5.1}
\begin{cases}
\square w
+\sum_{j,k=0}^{n+d}\gamma^{jk}(t,x,y)\partial_j\partial_k w + m^2
w=F, \quad 2B\le t\le T
\\
\partial_\nu w(t,x,y)=0, \quad y\in \partial\Omega
\\
w(t,x,y)=0, \quad t\le 2B.
\end{cases}
\end{equation}
The necessary assumption is that
\begin{multline}\label{5.2}
\sum_{0\le j,k\le n+d}\gamma^{jk}(t,x,y)\xi_j\theta_k=0, \, \,
\text{if } \, y\in \partial\Omega,
\\
\theta=(0,\dots,0,\nu_1(y),\dots\nu_d(y)), \, \text{and } \,
%\langle \, \xi,\nu(y)\, \rangle =0.
\xi\cdot\nu(y)=0.
\end{multline}
Under this assumption, the proof of Lemma 7.4.1 in \cite{H} shows
that \eqref{4.3} is valid when $w$ is a solution of \eqref{5.1}.
The assumption \eqref{1.8'} for the nonlinear terms implies that
if $\partial_\nu u(t,x,y)=0$, $y\in \partial\Omega$ then
$$\gamma^{jk}(t,x,y)=\sum_{l=0}^{n+d}A^{jk}_l\partial_l u +
uA^{jk}$$ must satisfy \eqref{5.2}.

Since we have just seen that we can apply \eqref{4.3}, \eqref{4.5} under the
assumption \eqref{1.8'}, we conclude that standard arguments give
Theorem \ref{theorem1.2}, i.e., the existence results for our
nonlinear Klein-Gordon equations with Neumann boundary conditions.

For Theorem \ref{theorem1.3}, which concerns existence for
semilinear wave equations with Neumann boundary conditions, one
just needs to use the special case of \eqref{4.3}, \eqref{4.5} where all of the
$\gamma^{jk}$ vanish identically, i.e., the standard energy
estimate for constant coefficients and Neumann conditions.  By
using this and Proposition \ref{prop3.6} one obtains Theorem
\ref{theorem1.3} using arguments from \cite{KSS2}.  The remainder of this section
will be dedicated to sketching this proof.

\noindent{\bf Proof of Theorem \ref{theorem1.3}:}  Using the reduction outline at the beginning
of the preceding section, it suffices to show that the solution $w$ to
\begin{equation}\label{5.3}
\begin{cases}
\Box w = (1-\eta) Q(\nabla_x(u_0+w))-[\Box,\eta](u_0+w)\\
\partial_\nu w(t,x,y)=0,\quad y\in\partial\Omega\\
w(t,x,y)=0,\quad t\le 2B
\end{cases}
\end{equation}
exists for $t\in [0,T_\varepsilon)$ when $n=3$ and exists globally for $n\ge 4$.
We solve this equation using iteration.
%Let us focus on the $n=3$ case.
To do so, we then let $w_0 \equiv 0$ and set $w_k$ to be the solution of
\begin{equation}\label{5.4}
\begin{cases}
\Box w_k = (1-\eta)Q(\nabla_x(u_0+w_{k-1}))-[\Box,\eta](u_0+w_k)\\
\partial_\nu w_k(t,x,y)=0,\quad y\in\partial\Omega\\
w_k(t,x,y)=0,\quad t\le 2B
\end{cases}
\end{equation}
for $k=1,2,3,\dots$.

We will focus on providing a sketch of the proof of Theorem \ref{theorem1.3} for the $n=3$
case using arguments of \cite{KSS2}.
Similar iteration techniques can be used to prove the $n\ge 4$; see, e.g., \cite{Metcalfe}.

We let
\begin{multline}\label{5.5}
M_k(T)=\sup_{0\le t\le T} \sum_{|\alpha|\le 10} \Bigl(\|\Gamma^\alpha \nabla_{t,x} w_k(t,\cd)\|_{L^2(\R^3\times\Omega)}
\\+ (\ln(2+t))^{-1/2} \|\langle x\rangle^{-1/2} \Gamma^\alpha \nabla_x w_k\|_{L^2([0,t]\times\R^3\times\Omega)}\Bigr).
\end{multline}
Our first goal is to inductively prove that for $\varepsilon$ sufficiently small we have
\begin{equation}\label{5.6}
M_k(T)\le 4C_0\varepsilon
\end{equation}
for a uniform constant $C_0$ greater than the constant occuring in \eqref{4.2},
$T\in [2B,T_\varepsilon)$, and $k=1,2,3,\dots$.  When $k=1$, \eqref{5.6} follows
from \eqref{1.6}, well-known local estimates, \eqref{3.10}, and Gronwall's inequality.

Let us prove \eqref{5.6} inductively.  Thus, we assume that this bound holds for $k-1$.
By \eqref{3.10} and \eqref{4.2},
\begin{multline}\label{5.7}
M_k(T_\varepsilon) \le C \sum_{|\alpha|\le 10} \int_0^{T_\varepsilon}
\|\Gamma^\alpha Q(\nabla_x(u_0+w_{k-1}))(t,\cd)\|_{L^2(\R^3\times\Omega)}\:dt\\
+C\sum_{|\alpha|\le 9} \Bigl(\|\Gamma^\alpha
Q(\nabla_x(u_0+w_{k-1})(T_\varepsilon,\cd)\|_{L^2(\R^3\times\Omega)}
\\+(\log(2+T_\varepsilon))^{-1/2}\|\langle x\rangle^{-1/2} \Gamma^\alpha Q(\nabla_x(u_0+w_{k-1}))
\|_{L^2([0,T_\varepsilon]\times\R^3\times
\Omega)}\Bigr)\\
2C_0\varepsilon + C\int_{2B+(1/2)}^{2B+1} \|\Gamma^\alpha \partial_t w_k(t,\cd)\|_{L^2(\R^3\times\Omega)}\:dt.
\end{multline}

Since $Q$ is quadratic, we can apply (2.9) and the standard
Sobolev lemma to see that
\begin{multline*}
\sum_{|\alpha|\le 10} \|\Gamma^\alpha
Q(\nabla_x(u_0+w_{k-1})(t,\cd)\|_{L^2(\{2^j\le |x|\le
2^{j+1}\}\times \Omega)}
\\\le
C \sum_{|\alpha|\le 10} \|\langle x\rangle^{-1/2} \Gamma^\alpha
\nabla_x u_0(t,\cd)\|^2_{L^2(\{2^j\le |x| \le
2^{j+1}\}\times\Omega)}
\\+C\Bigl(\sum_{|\alpha|\le 10} \|\Gamma^\alpha \nabla_x u_0(t,\cd)\|_{L^2(\R^3\times\Omega)}\Bigr)
\Bigl(\sum_{|\alpha|\le 10} \|\Gamma^\alpha \nabla_x
w_{k-1}(t,\cd)\|_{L^2(\{2^j\le |x|\le 2^{j+1}\}\times\Omega)}
\Bigr)
\\+C\Bigl(\sum_{|\alpha|\le 10} \|\Gamma^\alpha \nabla_x u_0(t,\cd)\|_{L^2(\{2^j\le |x|\le 2^{j+1}\}\times \Omega)}\Bigr)
\Bigl(\sum_{|\alpha|\le 10} \|\Gamma^\alpha \nabla_x
w_{k-1}(t,\cd)\|_{L^2(\R^3\times\Omega)}\Bigr)
\\+C \sum_{|\alpha|\le 10} \|\langle x\rangle^{-1/2} \Gamma^\alpha \nabla_xw_{k-1}(t,\cd)\|^2_{L^2
(\{2^{j-1}\le |x|\le 2^{j+2}\}\times\Omega)}.
\end{multline*}
Thus, since $u_0$ vanishes for $t\ge 2B+1$, summing over $j$ and applying \eqref{4.2} gives
$$\int_0^{T_\varepsilon} \|\Gamma^\alpha Q(\nabla_x(u_0+w_{k-1})(t,\cd)\|_{L^2(\R^3\times\Omega)}\:dt
\le C[\varepsilon+(\ln(2+T_\varepsilon))^{1/2} M_{k-1}(T_\varepsilon)]^2.$$
Using similar techniques, one can show that the second and third terms satisfy better bounds and are controlled by
$C(\varepsilon + M_{k-1}(T_\varepsilon))^2$.  Hence, by Gronwall's inequality and the inductive hypothesis \eqref{5.6}, we have
$$M_{k}(T_\varepsilon)\le 2C_0\varepsilon + 3C[\varepsilon+(\ln(2+T_\varepsilon))^{1/2}4C_0\varepsilon]^2.$$
Thus, one can choose an appropriate $c$ in \eqref{1.8} so that this is bounded by $4C_0\varepsilon$ as desired.

To show that the sequence $w_k$ converges to a solution, we estimate
\begin{multline}\label{5.8}
A_k(T)=\sup_{0\le t\le T} \sum_{|\alpha|\le 10} \Bigl(\|\Gamma^\alpha \nabla_{t,x}(w_k - w_{k-1})(t,\cd)\|
_{L^2(\R^3\times\Omega)}
\\+(\ln(2+t))^{-1/2}\|\langle x\rangle^{-1/2}
\Gamma^\alpha
\nabla_x(w_k-w_{k-1})\|_{L^2([0,t]\times\R^3\times\Omega)}\Bigr).
\end{multline}
If one shows that $A_k$ satisfies
\begin{equation}\label{5.9}
A_k(T_\varepsilon)\le \frac{1}{2}A_{k-1}(T_\varepsilon),\quad k=1,2,3,\dots,\end{equation}
the proof of Theorem \ref{theorem1.3} would be complete.  Indeed this is the case.  Since $Q$ is
quadratic, we can repeat the previous arguments to see that
$$A_k(T_\varepsilon)\le C \Bigl(2C_0\varepsilon+(\ln(2+T_\varepsilon))^{1/2}(M_{k-1}(T_\varepsilon)
+M_{k-2}(T_\varepsilon))\Bigr) (\ln(2+T_\varepsilon))^{1/2} A_{k-1}(T_\varepsilon).$$
By our choise of $T_{\varepsilon}$, we see that this leads to \eqref{5.9} for $\varepsilon$ sufficiently
small, which completes the proof of the $n=3$ case.

In order to show global existence for $n\ge 4$, we instead use \eqref{3.12} to bound
\begin{equation}\label{5.10}
M_k(T)=\sup_{0\le t\le T} \sum_{|\alpha|\le n+10}
\Bigl(\|\Gamma^\alpha \nabla_{t,x} w_k(t,\cd)\|_2 +\|\langle
x\rangle^{-(n-1)/4} \Gamma^\alpha
\nabla_xw_k\|_{L^2([0,t]\times\R^n\times\Omega)}\Bigr)
\end{equation}
for any $T\in [2B,\infty)$.  Similarly, we show that
\begin{multline}\label{5.11}
A_k(T)=\sup_{0\le t\le T} \sum_{|\alpha|\le n+10}
\Bigl(\|\Gamma^\alpha \nabla_{t,x}
(w_k-w_{k-1})(t,\cd)\|_{L^2(\R^n \times\Omega)} \\+ \|\langle
x\rangle^{-(n-1)/4} \Gamma^\alpha \nabla_x
(w_k-w_{k-1})\|_{L^2([0,t]\times\R^n\times \Omega)}\Bigr)
\end{multline}
is a Cauchy sequence, which would complete the proof.\qed

\end{document}